\documentclass[11pt]{article}
\usepackage{amsmath,amssymb}
\usepackage{epsfig}
\usepackage{amsfonts}
\usepackage{amsmath,amssymb}
\usepackage{epsfig}
\def\dref#1{(\ref{#1})}

\begin{document}

\vspace*{1\baselineskip} \begin{center}\Large\bf Analysis and
control of network synchronizability \footnote {\small
 This work is supported by the
 National
 Science Foundation of
 China under grants 60674093, 60334030.}
 \end{center}
\vspace*{1\baselineskip}

\centerline {Zhisheng Duan,  ~Guanrong Chen ~and ~Lin Huang}

\vspace*{0.5\baselineskip}
\begin{center}
 State Key Laboratory for Turbulence and Complex Systems,
Department of Mechanics and Aerospace Engineering, College of
Engineering, Peking University, Beijing 100871, P. R. China
 \\ Emails:
       duanzs@pku.edu.cn,  eegchen@cityu.edu.hk, hl35hj75@pku.edu.cn  \\
{\it Tel/Fax}: (8610)62765037
\end{center}
\vskip 0.5cm

{\bf Abstract.} \,\, In this paper, the investigation is first
motivated by showing two examples of simple regular symmetrical
graphs, which have the same structural parameters, such as average
distance, degree distribution and node betweenness centrality, but
have very different synchronizabilities. This demonstrates the
complexity of the network synchronizability problem. For a given
network with identical node dynamics, it is further shown that two
key factors influencing the network synchronizability are the
network inner linking matrix and the eigenvalues of the network
topological matrix. Several examples are then provided to show that
adding new edges to a network can either increase or decrease the
network synchronizability. In searching for conditions under which
the network synchronizability may be increased by adding edges, it
is found that for networks with disconnected complementary graphs,
adding edges never decreases their synchronizability. This implies
that better understanding and careful manipulation of the
complementary graphs are important and useful for enhancing the
network synchronizability. Moreover, it is found that an unbounded
synchronized region is always easier to analyze than a bounded
synchronized region. Therefore, to effectively enhance the network
synchronizability, a design method is finally presented for the
inner linking matrix of rank 1 such that the resultant network has
an unbounded synchronized region, for the case where the synchronous
state is an equilibrium point of the network.

{\bf Keywords.} \,\, Complex network, Complementary graph,
Synchronizability, Synchronized region, Inner linking matrix.

\section{Introduction and problem formulation}

 The subject of network synchronization has recently attracted increasing attention
 from various fields
(see \cite{bar02,bel05, boc06, cur97,  luw04, lu04,  sor07, wang06,
wang02, wat98, wu02, zhao06}  and references therein).
 Of particular importance is how the synchronizability depends on
 various structural parameters of the network, such as average
  distance,
 clustering coefficient, coupling strength, degree distribution and weight distribution,
 among others.
Some important results have been established for such concerned
problems based on the notions of master stability function and
synchronized region \cite{bar02,koc05,lu04,mot05,pec98,zhou06}. Some
interesting relationships between  synchronizability and structural
parameters of networks have also been reported, e.g., smaller
average network distance does not necessarily mean better
synchronizability \cite{nis03}, therefore the betweenness centrality
was proposed as a good indicator for synchronizability
\cite{hong04}. And two networks with the same degree sequence were
constructed in a probabilistic sense to demonstrate that they can
have different synchronizabilities  \cite{wu05}, showing that
synchronizability has no direct relations with degree distribution.
Moreover, the effect of  perturbations of coupling matrices on the
synchronizability was studied in \cite{wu03}.
 Motivated by all these research works, this paper attempts to further explore the analysis and control problems
 of synchronizability for various complex dynamical
networks.

Consider a dynamical network consisting of $N$ coupled identical
nodes, with each node being an $n$-dimensional dynamical system,
described by
\begin{equation}\label{n1}
\dot{x}_i=f(x_i)-c\sum_{j=1}^Na_{ij}H(x_j),\;i=1,2,\cdots,N,
\end{equation}
where $x_i=(x_{i1},x_{i2},\cdots,x_{in})\in \mathbb{R}^n$ is the
state vector of node $i$, $f(\cdot):\mathbb{R}^n\rightarrow
\mathbb{R}^n$
is a smooth 
vector-valued function, constant $c>0$ represents the
\textit{coupling strength},  $H(\cdot):\mathbb{R}^n\rightarrow
\mathbb{R}^n$ is called the
\textit{inner linking function}, 
and $A=(a_{ij})_{N\times N}$ is called the \textit{outer coupling
matrix} or \textit{topological  matrix}, which represents the
coupling configuration of the entire network. This paper only
considers the case that  the network is diffusively connected, i.e.,
$A$ is irreducible and its entries satisfy
$$
a_{ii}=-\sum_{j=1,j\neq i}^Na_{ij},\;i=1,2,\cdots,N.
$$
Further, suppose that, if there is an edge between node $i$ and node
$j$, then $a_{ij}=a_{ji}=-1$, i.e., $A$ is a Laplacian matrix.
Therefore, $0$ is an eigenvalue of $A$ with multiplicity 1, and all
the other eigenvalues of $A$ are strictly positive, which are
denoted by
 \begin{equation}\label{f1}
0=\lambda_1<\lambda_2\leq\lambda_3\leq\cdots\leq\lambda_N.
\end{equation}

The dynamical network \dref{n1} is said to achieve (asymptotical)
synchronization if
\begin{equation}\label{f2}
x_1(t)\rightarrow x_2(t)\rightarrow\cdots\rightarrow
x_N(t)\rightarrow s(t),\;\textrm{as}\; t\rightarrow \infty,
\end{equation}
where, because of the diffusive coupling configuration,  the
\textit{synchronous state} $s(t)\in \mathbb{R}^n$ is a solution of
an individual node, i.e., $\dot{s}(t)=f(s(t))$. Here, $s(t)$ can be
an equilibrium point, a periodic orbit, or even a chaotic orbit.

As shown in \cite{lu04,pec98}, the local stability of the
synchronized solution $x_1(t)=x_2(t)=\cdots=x_N(t)=s(t)$
can be determined by analyzing the following so-called
\textit{master stability equation}:
\begin{equation}\label{f3}
\dot{\omega}=[Df(s(t))+\alpha DH(s(t))]\omega,
\end{equation}
where $\alpha\in \mathbb{R}$, and $Df(s(t))$ and $DH(s(t))$ are the
Jacobian matrices of functions $f$ and $H$ at $s(t)$, respectively.

The largest Lyapunov exponent $L_{max}$ of network \dref{n1}, which
can be calculated from system \dref{f3} and is a function of
$\alpha$, is referred to as the \textit{master stability function}.
In addition, the region $S$ of negative real $\alpha$ where
$L_{max}$ is also negative is called the \textit{synchronized
region} of network \dref{n1}. Based on the results of
\cite{lu04,pec98},
 the synchronized solution
 of dynamical network \dref{n1} is locally asymptotically stable if, and
only if,
\begin{equation}\label{f4}
c\lambda_k\in S,\;k=2,3,\cdots,N.
\end{equation}
 The synchronized region $S$ can be an
unbounded region, a bounded region,  an empty set, or a union of
several such regions.

Obviously, for given node dynamics of a linearly coupled network,
two key factors influencing the synchronizability are the inner
linking matrix $H(\cdot)=H$ and the eigenvalues of the topological
matrix $A$. The inner linking matrix is directly related to the
synchronized region, as studied in \cite{duan07, koc05, liu07,
pec98, ste07}. The larger the synchronized region, the easier the
synchronization. The topological  matrix, on the other hand, is
directly related to condition (\ref{f4}). If $S$ is an unbounded
sector $(-\infty, \alpha]$, the eigenvalue $\lambda_2$ of $A$
determines the synchronizability \cite{wang02}; if $S$ is a bounded
sector $[\alpha_1, \alpha_2]$, the ratio
$r(A)=\frac{\lambda_2}{\lambda_N}$ determines the synchronizability
\cite{bar02}. No matter what the synchronized region is, the larger
the $\lambda_2$ and $r(A)$ are, the easier the synchronization is.
This paper will further study this issue with more careful analysis.

The rest of this paper is organized as follows. In Section 2, two
simple graphs on six nodes are given to show  that  networks with
the same structural parameters, such as average
 distance,  degree distribution and betweenness centrality, can have different
 synchronizabilities.
   If the synchronized region $S$ is unbounded, adding edges never
decreases the synchronizability, but this may not be true if $S$ is
bounded. In Section 3, a class of networks with disconnected
complementary graphs are discussed. For such networks, adding edges
never decreases the synchronizability no matter what type of region
$S$ is. In Section 4, an inner linking matrix of rank 1 is designed
for realizing unbounded synchronized regions in the case that the
synchronous state is an equilibrium point. In Section 5, some
network synchronization examples are provided to illustrate the
theoretical results. The paper is concluded by the last section.

Throughout this paper, for any given undirected graph $G$,
eigenvalues of $G$ mean eigenvalues of its corresponding Laplacian
matrix. Notations for graphs and their corresponding Laplacian
matrices are not differentiated, and networks and their
corresponding graphs are not distinguished, unless otherwise
indicated.

\section{Two simple graphs that tell the main idea}

In \cite{wu05}, two networks with the same degree sequence were
constructed in a probabilistic sense  to show that they can have
different synchronizabilities. In this section,  the  two simple
graphs $G_1$ and $G_2$ on six nodes, shown in Figs. 1 and 2, are
considered, where $G_1$ is a typical bipartite graph with many
interesting properties.

\begin{center}
 \unitlength=1cm
 \qquad \hbox{\hspace*{0.1cm}  \epsfxsize5.5cm \epsfysize5cm
\epsffile{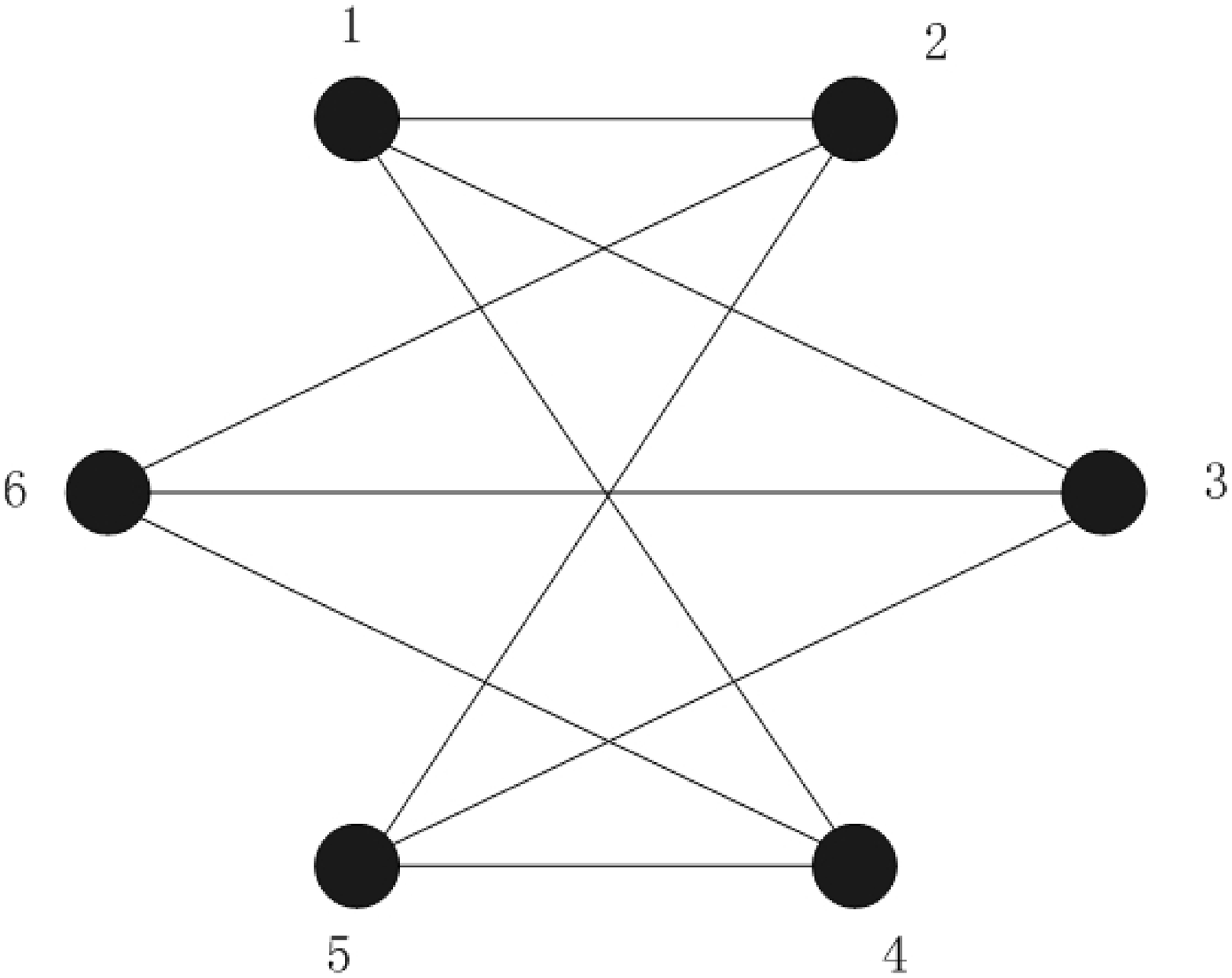}\qquad \epsfxsize5.5cm \epsfysize5cm
\epsffile{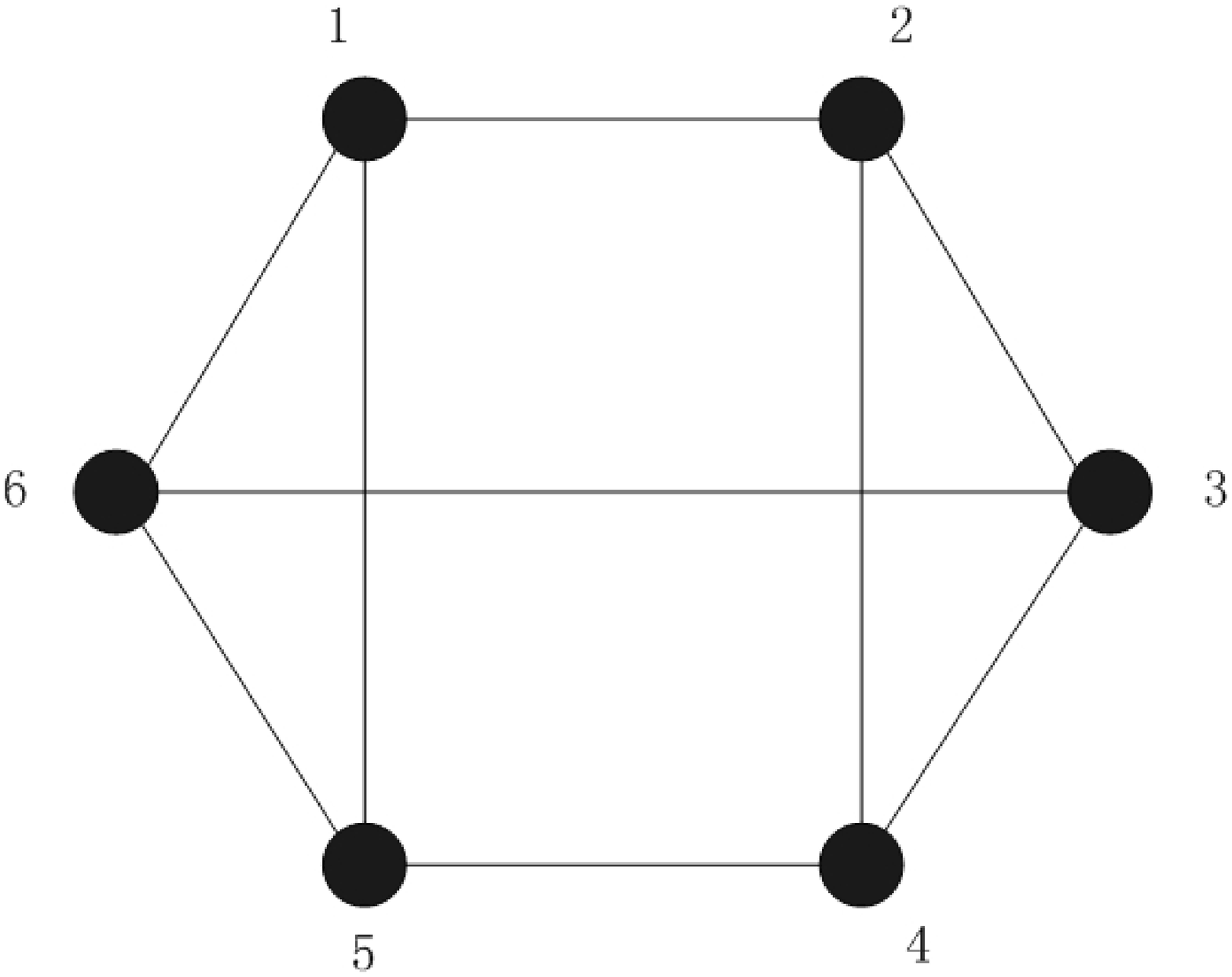}
 }
\end{center}
\vskip -0.3cm \centerline{  Fig. 1 \,\, Graph $G_1$
\qquad\qquad\qquad \qquad Fig. 2 \,\, Graph $G_2$}

Obviously, graphs $G_1$ and $G_2$ have the same degree sequence,
where the degree of every node is 3; the same average  distance
$\frac{7}{5}$; and the same node betweenness centrality 2
\cite{hol02, hong04}.
 Although these two graphs have the same
 structural characteristics, their corresponding networks  have
 different synchronizabilities, as shown below.
Their Laplacian matrices  are
  $$G_1=\left(\begin{array}{cccccc}
  -3 &1 &1 & 1& 0& 0\\
    1 &-3& 0 &0& 1& 1\\
    1 &0& -3& 0& 1& 1\\
    1 &0& 0 &-3& 1 &1\\
    0 &1& 1& 1& -3& 0\\
    0 &1 &1 &1& 0& -3 \end{array}\right),\quad
  G_{2}=\left(\begin{array}{cccccc}
  -3 &1& 0& 0& 1& 1\\
    1& -3& 1 &1& 0& 0\\
    0 &1& -3& 1& 0 &1\\
    0 &1 &1& -3& 1& 0\\
    1 &0& 0 &1 &-3& 1\\
    1 &0& 1& 0& 1& -3 \end{array}\right),
$$
respectively. The eigenvalues of $G_1$ are $0, 3, 3, 3, 3 $ and $6$;
the eigenvalues of $G_2$ are $0, 2, 3, 3, 5 $ and $5$. Obviously,
$\lambda_2(G_1)=3>\lambda_2(G_2)=2,$ and $r(G_1)=0.5>r(G_2)=0.4$.
Therefore,  the synchronizability of network $G_1$ is better than
that of network $G_2$.

 Graphs $G_1$ and $G_2$ have the same structural parameters,
  but it is clear that they have different
average clustering  coefficients, denoted by $C(G_i), i=1,2,$
 with $C(G_2)>C(G_1)$. As mentioned above,
 the clustering coefficient  does not have direct relation to synchronization \cite{nis03}.
 For example,
 globally coupled graphs have the largest clustering coefficient, 1,
 and they have the best synchronizability. However, for the
 above two graphs, the larger clustering coefficient does not indicate
  better synchronizability. This is  demonstrated by
 the following process of adding edges.

 Consider enhancing $\lambda_2$ and $r$
 by adding edges to $G_2$.  For this purpose, the following lemma is
 needed.

 {\bf Lemma 1}\,\cite{mer98}\, For any given connected
 undirected graph $G$ of size $N$, its nonzero
  eigenvalues indexed as in
  (\ref{f1}) grow monotonically with the
 number of added edges, that is, for any added edge $e$,
 $\lambda_i(G+e)\geq \lambda_i(G)$, $i=1, \cdots, N$.

By Lemma 1, obviously if the synchronized region is unbounded,
adding edges never decreases the synchronizability. However, for
bounded synchronized regions, this is not necessarily true. For
example, adding an edge between node 1 and node 3 in
 graph $G_2$ (Fig. 2), denoted by $e\{1,3\}$, leads to a new graph
 $G_2+e\{1,3\}$,
 whose eigenvalues are $0, 2.2679, 3, 4, 5$ and $5.7321$.
Thus, $ r(G_2+e\{1,3\})=0.3956$ is even smaller than $r(G_2)=0.4$.
This means that the synchronizability of network $G_2+e\{1,3\}$ is
worse than that of network $G_2$. Adding a new edge between node 1
and node 4 instead, one gets $
r(G_2+e\{1,3\})<r(G_2+e\{1,3\}+e\{1,4\})=0.3970<r(G_2).$ This means
that, the synchronizability of network $G_2+e\{1,3\}+e\{1,4\}$ is
better than $G_2+e\{1,3\}$, but still worse than $G_2$. Therefore,
by adding edges, the network synchronizability may increase or
decrease, for which no general rule has been found to date.

On the other hand,  during the process of adding edges, average
distance decreases and average clustering coefficient increases. But
this does not indicate better synchronizability, consistent with the
conclusion in \cite{nis03}.

It was shown \cite{hong04} that the synchronizability is always
improved as the maximum betweenness centrality is reduced, which is
consistent with the conclusion of \cite{nis03}. In the above two
graphs, however, it shows that the same betweenness centrality does
not necessarily mean the same synchronizability. On the other hand,
adding three edges between nodes 1 and 6, 2 and 3, 3 and 4,
respectively,  in graph $G_1$, and then computing their
corresponding  eigenvalues, it can be verified that the networks
corresponding to the resulting graphs
$G=G_1+e\{1,6\}+e\{2,3\}+e\{3,4\}$ and $G_1$ have the same
synchronizability. However, in this case, the maximum betweenness
centrality of $G$, $\frac{11}{6}$, is smaller than that of $G_1$,
$2$. This shows that the smaller betweenness centrality does not
necessarily indicate better synchronizability, revealing the
complexity in the relationship between synchronizability and network
structural parameters.

Note that adding edges in $G_1$ also increases the clustering
coefficient and decreases the average distance, but this does not
result in the increase of synchronizability. In the following, it
explains why adding three edges in $G_1$ does not increase the
synchronizability.

\section{Networks with disconnected complementary graphs}

For a given graph $G$, the complement of $G$ is the graph containing
all the nodes of $G$, and all the edges that are not in $G$. The
complementary graph of $G$ is  denoted by $G^c$. For example, the
complementary graphs  of $G_1$ and $G_2$ in Fig. 1 and Fig. 2 are
shown in Fig. 3 and Fig. 4, respectively. In the previous section,
it shows that adding edges sometimes decreases the
synchronizability. However, for a class of graphs with disconnected
complementary graphs, this never occurs. In order to discuss such
networks, the following lemma is needed.

\begin{center}
 \unitlength=1cm
 \qquad \hbox{\hspace*{0.1cm}  \epsfxsize5.5cm \epsfysize5cm
\epsffile{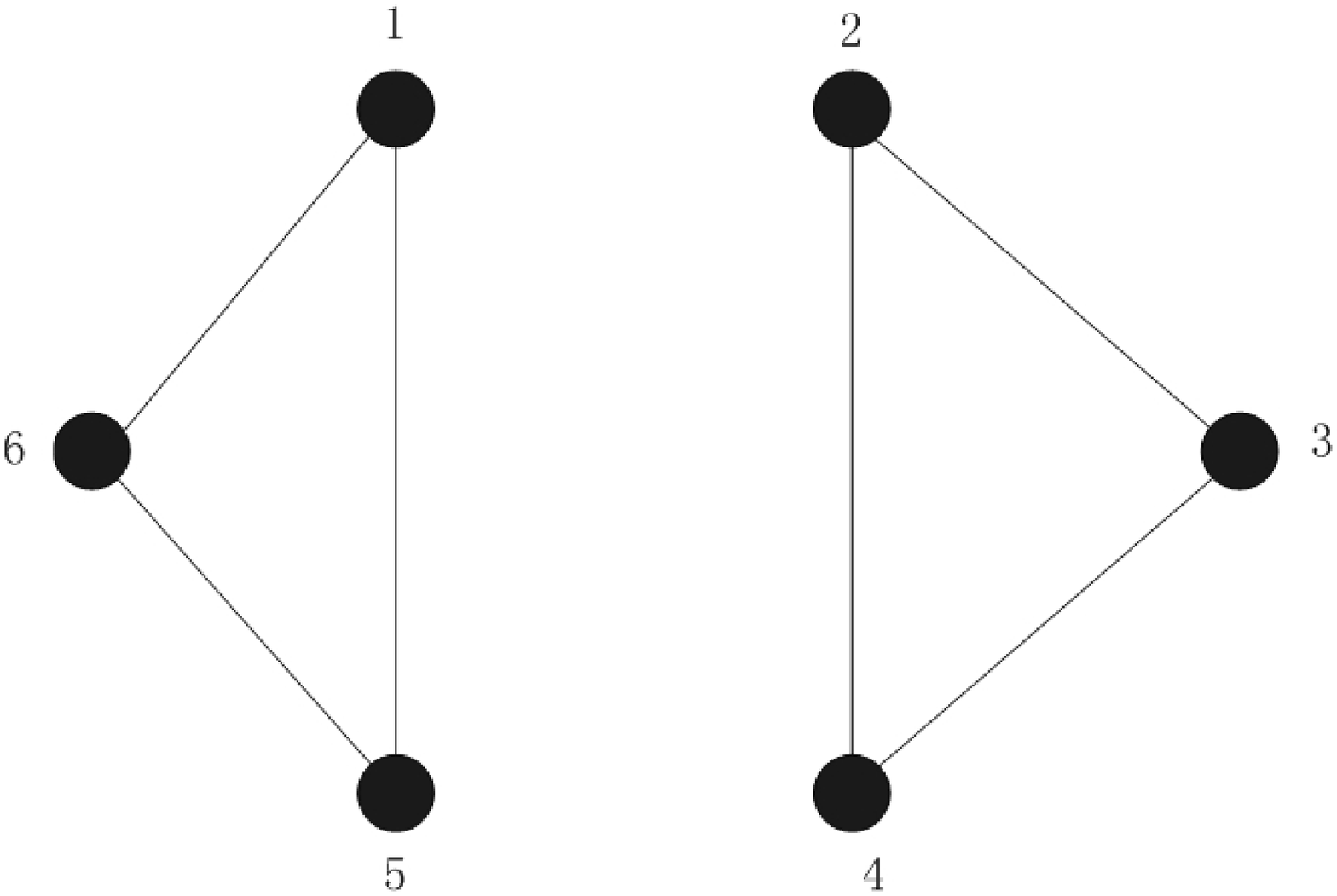}\qquad \epsfxsize5.5cm \epsfysize5cm
\epsffile{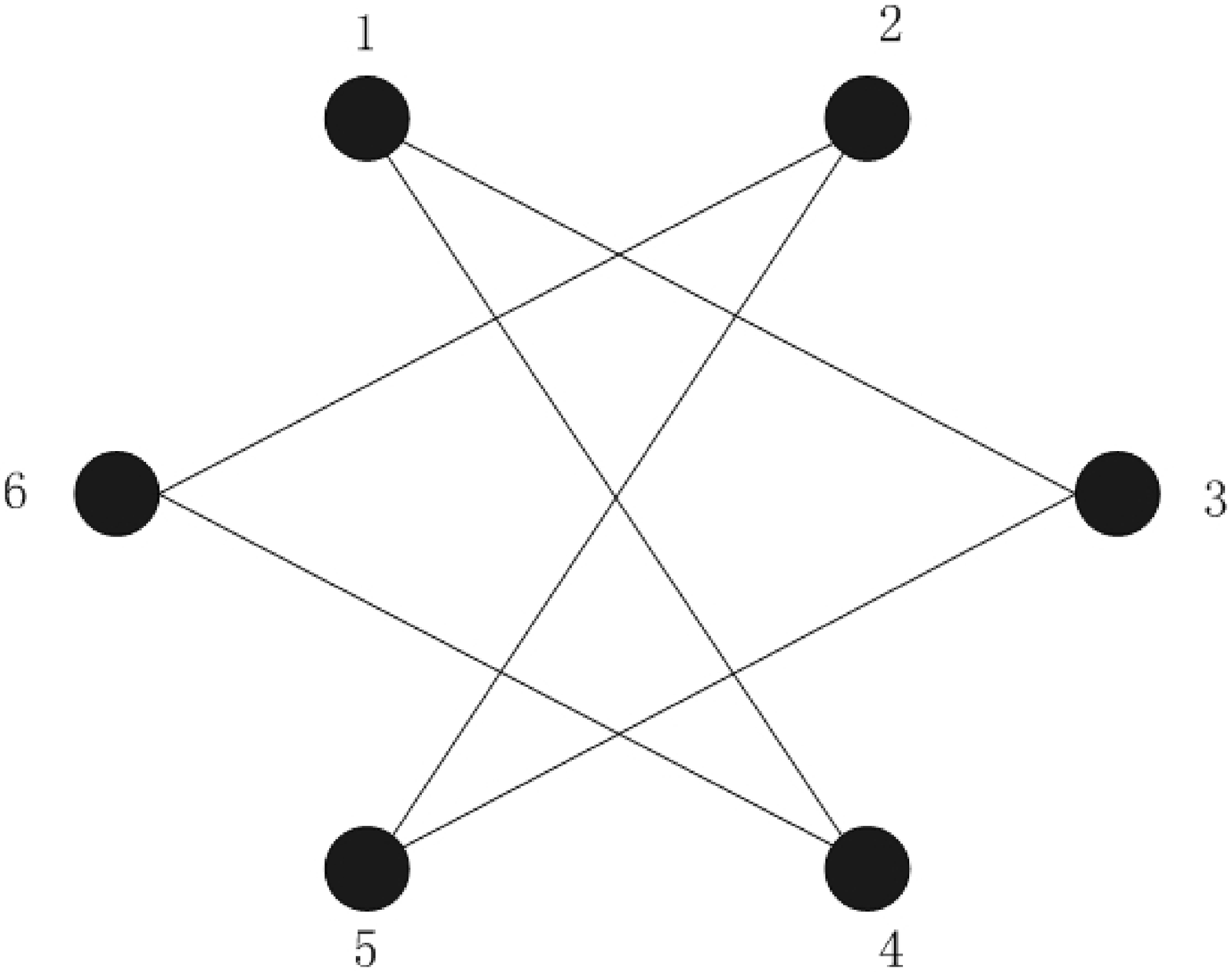}
 }
\end{center}
\vskip -0.3cm \centerline{  Fig. 3 \,\, Graph $G_1^c$
\qquad\qquad\qquad \qquad Fig. 4 \,\, Graph $G_2^c$}

 {\bf Lemma 2}\,\cite{mer98}\, For any given graph $G$, the following statements
 hold:

(i) \, $\lambda_N(G)$, the largest eigenvalue of $G$, satisfies
$\lambda_N(G)\leq N.$

(ii) \, $\lambda_N(G)= N$ if, and only if, $G^c$ is disconnected.

(iii)\, If $G^c$ is disconnected and  has (exactly) $q$ connected
components, then the multiplicity of $\lambda_N(G)=N$  is $q-1$.

(iv) \, $\lambda_i(G^c)=N-\lambda_{N-i+2}(G),\quad 2\leq i\leq N$.

The complementary graph of $G_1$ (Fig.1) is shown in Fig. 3, which
is disconnected. The largest eigenvalue of $G_1$ is 6, which remains
the same when the graph receives additional  edges. Hence, combining
with Lemma 1, the synchronizability of the networks built on graph
$G_1$ never decrease with adding edges. Although this is true,
adding any 3 edges to graph $G_1$ does not enhance the
synchronizability, since the least nonzero eigenvalue $\lambda_2=3$
of $G_1$ has multiplicity 4 (the multiplicity of the largest
eigenvalue in $G_1^c$). This is due to the fact that, for any graph
$G$, $\hbox{rank} (\lambda_iI-(G+e))\leq  \hbox{rank}
(\lambda_iI-G)+1$.

According to Lemma 2, the multiplicity of the largest eigenvalue of
a graph $G$ is related to the number of connected components of its
complement $G^c$. Hence, with the same number of edges, generally
the synchronizability of  the networks built on graph $G$ is better
when $G^c$ has two connected components than the case that $G^c$ has
more than two connected components. This is because in the latter
case, the multiplicity of the largest eigenvalue of $G$ is larger
than 1, i.e., some edges have no contributions to the
synchronizability. On the other hand, the multiplicity of the
largest eigenvalue of $G^c$ (i.e., the multiplicity of the least
nonzero eigenvalue of $G$) should be large in order to reduce  the
number of edges needed to enhance the synchronizability. Therefore,
better understanding and careful manipulation of complementary
graphs are useful for enhancing the network synchronizability. And,
at least for dense networks, the complementary graphs are easier to
analyze than the original graphs, e.g., $G_1^c$ in Fig. 3 is simpler
than $G_1$ in Fig. 1.

{\bf Remark 1}\,\, The graphs shown in Figs. 1 and 2 can be
generalized to graphs of size $N=2n$. Suppose that graph $G_1$ is
bipartite in the sense that it contains two sets of nodes, each set
containing $n$ isolated nodes, and each node in one set connects to
all the nodes in the other set, i.e., the complementary graph of
$G_1$ is two separated fully connected subgraphs of size $n$. Graph
$G_2$ is composed of two fully connected subgraphs and $n$ edges
connecting each node in one subgraph to the corresponding  node in
the other subgraph. In this case, the least nonzero and maximum
eigenvalues of $G_1$ are $\frac{N}{2}$ and $N$, respectively, and
$r(G_1)=\frac{1}{2}.$ On the other hand, the least nonzero and
maximum eigenvalues of $G_2$ are $2$ and $\frac{N}{2}+2$
\cite{mer98}, respectively, with $r(G_2)=\frac{4}{N+4}\rightarrow 0$
as $N\rightarrow +\infty$. Therefore, these two graphs have  the
same structural parameters but have very different
synchronizabilities.

\section{Designing the inner linking matrix}

From the above two sections, it can be seen that the bounded
synchronized regions are  more complicated than the unbounded
synchronized regions. Thus, the synchronizability is easier to
analyze when the synchronized region is unbounded. Hence, it is
interesting to find out  how to design the inner linking matrix such
that the network synchronized region is unbounded.

 If the synchronous state is
an equilibrium point, then both $Df(s(t))$ and $DH(s(t))$ in
(\ref{f3}) reduce to constant matrices, denoted by $F$ and $H$,
respectively. In this case, system \dref{f3} becomes
\begin{equation}\label{f5}
\dot{\omega}=[F+\alpha H]\omega.
\end{equation}
 Hence, the synchronized region $S$ becomes the stability region of
$F+\alpha H$ with respect to parameter $\alpha$.  In this section,
consider the design of an $H$ such that $F+\alpha H$ has an
unbounded stable region. It is well known that if $H$ is an
anti-stable matrix (e.g., $H=I_n$), $F+\alpha H$ has an unbounded
stable region \cite{liu07}. However, if $H$ is of full rank, it
means that the coupling in the network is a full state coupling
among nodes, so the coupling cost may be high. For this reason,
consider the design of an $H$ of rank 1 such that the stability
region for $F+\alpha H$ is unbounded. In this case, the coupling can
be viewed as an input-output coupling as in control systems
\cite{duan04}, or an observer-based coupling \cite{jia06}.

{\bf Theorem 1}\,\, Given a matrix $F\in {\bf R}^{n\times n}$, there
exists a matrix $H\in {\bf R}^{n\times n}$ of rank 1 such that the
stability region of $F+\alpha H$ with respect to parameter $\alpha$
contains $(-\infty, \alpha_1]$, $\alpha_1<0$, if and only if every
unstable eigenvalue of $F$ is corresponding to only one Jordan
block.

{\bf Proof }\,\, Without loss of generality, suppose $\alpha_1=-1.$
If $F+\alpha H$ is stable and $H$ is of rank 1, it means that
$(F,H)$ is stabilizable, so that every unstable eigenvalue of $F$
must be corresponding to only one Jordan block.

On the other hand,  if every unstable eigenvalue of $F$ is
corresponding to only one Jordan block, one may take a column vector
$b$ such that $(F, b)$ is stabilizable. Then, there exists a row
vector $k$ such that $F-bk$ is stable, i.e., there exists a matrix
$P=P^T>0$ such that
$$(F-bk)P+P(F-bk)'<0.$$
Let $kP=y$, so that the above inequality becomes
 \begin{equation}  \label{th1}
 FP+PF'-by-y'b'<0. \end{equation} By the
canonical projection lemma in $H_\infty$ control theory
\cite{iwa94}, there exists $y$ such that the above inequality holds
if, and only if, there exists a scalar $\beta>0$ such that
 \begin{equation}  \label{th2}
 FP+PF'-\beta bb'<0. \end{equation}
Since $P$ is a matrix to be determined, in the above inequality,
suppose $\beta=2$ without loss of generality. Then
\begin{equation}  \label{th3} FP+PF'-2 bb'<0.\end{equation}
Obviously, if $(F,b)$ is stabilizable, then there exists $P=P^T>0$
such that (\ref{th3}) holds. And, when (\ref{th3}) holds, for any
$\beta\geq 2$, (\ref{th2}) holds. Take $y=b'$, i.e., $k=b'P^{-1}$.
 By the above inequalities, one knows that $F-\alpha bk$ is stable
for all $\alpha\in (-\infty, -1].$ Therefore, $H=bk$ is  the matrix
to be found. \hfill $\Box$

{\bf Remark 2}\,\,Let $z_i=kx_i$ and the inner linking function
$H(x_j)=bkx_j$ in network (\ref{n1}). Then $z_i$ can be viewed as
the output of node $i$ of (\ref{n1}) and the linking function
$bkx_j$ can be viewed as the influence of the output of node $j$ to
the other nodes. Clearly, the above coupling is simpler than full
state couplings.

{\bf Remark 3}\,\, If $F$ is stable, i.e., the node system is
locally stable, then there is always a matrix $H$ of rank 1 such
that the resulting network has an unbounded synchronized region, as
shown by the examples given below. In addition, if there exists an
unstable eigenvalue of $F$ whose number of  Jordan blocks is more
than 1, one may take an inner coupling matrix $H$ such that its rank
equals the maximum number of Jordan blocks among the unstable
eigenvalues of $F$. On the other hand,  one may design an $H$ such
that $F+\alpha H$ has an unbounded stable region. Moreover, one may
also design an $H$ such that the unstable region of $F+\alpha H$ is
unbounded, if desired, which is useful for desynchronization
problems.

\section{Examples}

{\bf Example 1}\,\, Consider the network (\ref{n1}) consisting of
the third-order smooth Chua's circuits \cite{tsu05}, in which each
node  is described by
 \begin{equation}  \label {sc1}
\begin{array}{ccl}
 \dot{x}_{i1}&=&-k\alpha x_{i1}+k\alpha x_{i2}-k\alpha(a x_{i1}^3+bx_{i1}),\\
 \dot{x}_{i2}& =& kx_{i1}-kx_{i2}+kx_{i3},\\
 \dot{x}_{i3}& =& -k\beta x_{i2}-k\gamma x_{i3}. \end{array}\end{equation}
The vector $x_i$ in (\ref{n1}) is $(x_{i1}, x_{i2}, x_{i3})^T$ here.
Linearizing (\ref{sc1})  at its zero equilibrium gives
 \begin{equation}  \label {lsc1}
 \dot{x}_{i}= Fx_i, \quad F=\left(\begin{array}{ccc} -k\alpha-k\alpha b &k\alpha & 0\\
 k & -k & k\\0 & -k\beta &-k\gamma  \end{array}\right).
 \end{equation}

 Take $k=1,
\alpha=-0.1, \beta=-1, \gamma=1,  a=1, b=-25.$ Then $F$ is stable,
i.e., the node system (\ref{sc1}) is locally stable about zero.
Further, take the inner linking matrix \cite{duan07}
 $$H=\left(\begin{array}{ccc} 0.8348 &  9.6619&  2.6591\\ 0.1002 & 0.0694 & 0.1005\\
  -0.3254 &  -8.5837 &  -0.9042  \end{array}\right).$$
  Then, by simple computation, one knows that $F+\alpha H$ has two
  disconnected stable regions: $S_1=[-0.0099, 0]$ and $S_2=[-2.225, -1)$.
  Therefore, the entire synchronized region
  is $S_1\bigcup S_2$. Moreover, suppose that the number of nodes  is $N=6$,
and  the outer coupling matrix $A$ is equal
  to the $G_1$ in Section 2. According to the eigenvalues of $G_1$
  given in Section 2, one may take the coupling strength $c=\frac{1}{2.9}$.
  Then, for every eigenvalue of $G_1$, one has $c\lambda_i\in S_2$. By (\ref{f4}),
   network (\ref{n1}) specified with the above data achieves
   synchronization. However, for the outer coupling matrix $G_2$
  given in Section 2, for any coupling strength $c\in [0.002, +\infty)$,
   (\ref{f4}) does not hold. Therefore, for the above node equation, inner
  coupling matrix and coupling strength, the network built on graph $G_1$ in Fig. 1
  achieves synchronization, but the network built on $G_2$ in Fig. 2
  does not synchronize. Figs. 5 and 6 show their
 synchronous and non-synchronous behaviors, respectively.

\begin{center}
 \unitlength=1cm
 \qquad \hbox{\hspace*{0.1cm}  \epsfxsize5.5cm \epsfysize5cm
\epsffile{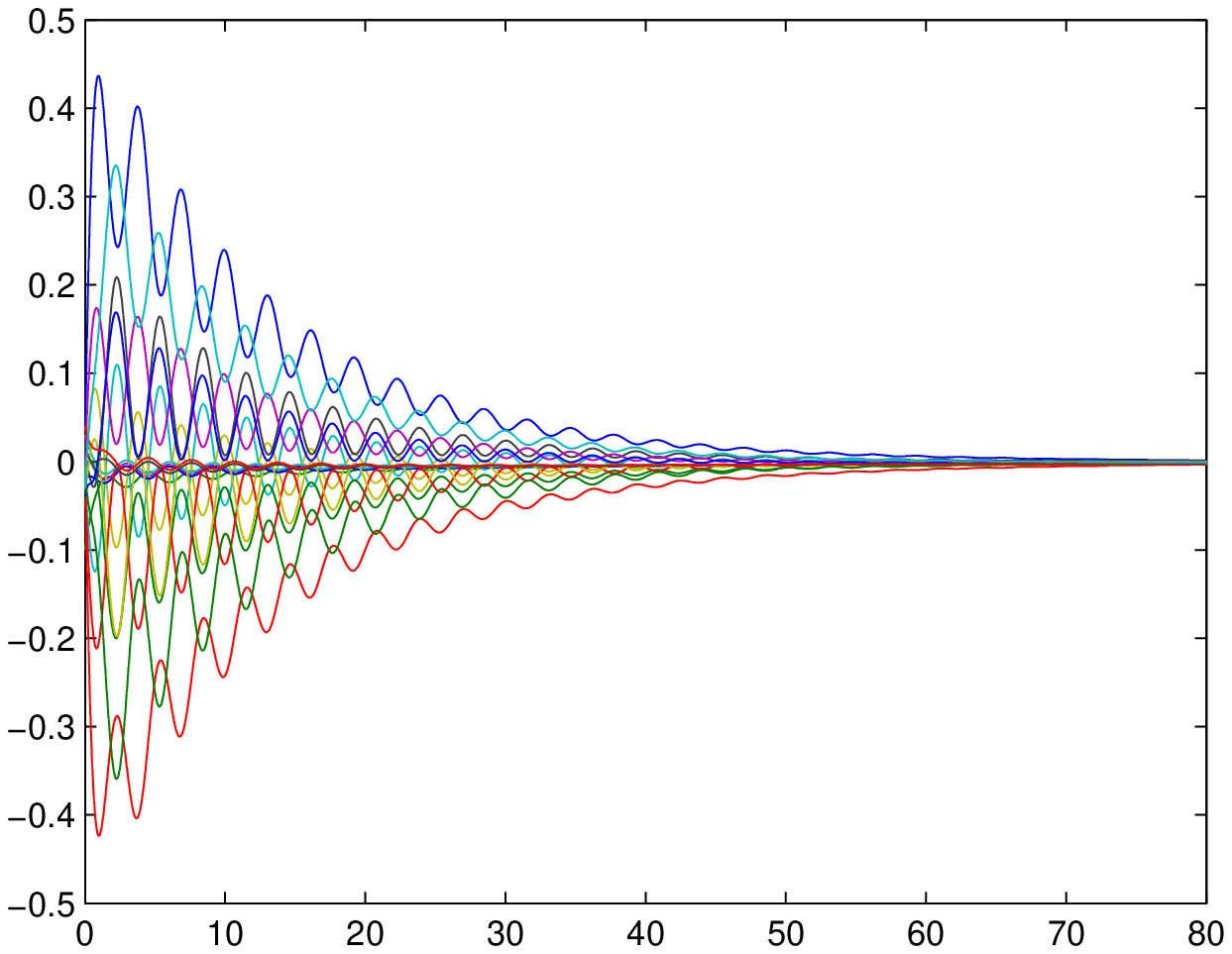}\qquad \epsfxsize6cm \epsfysize6cm
\epsffile{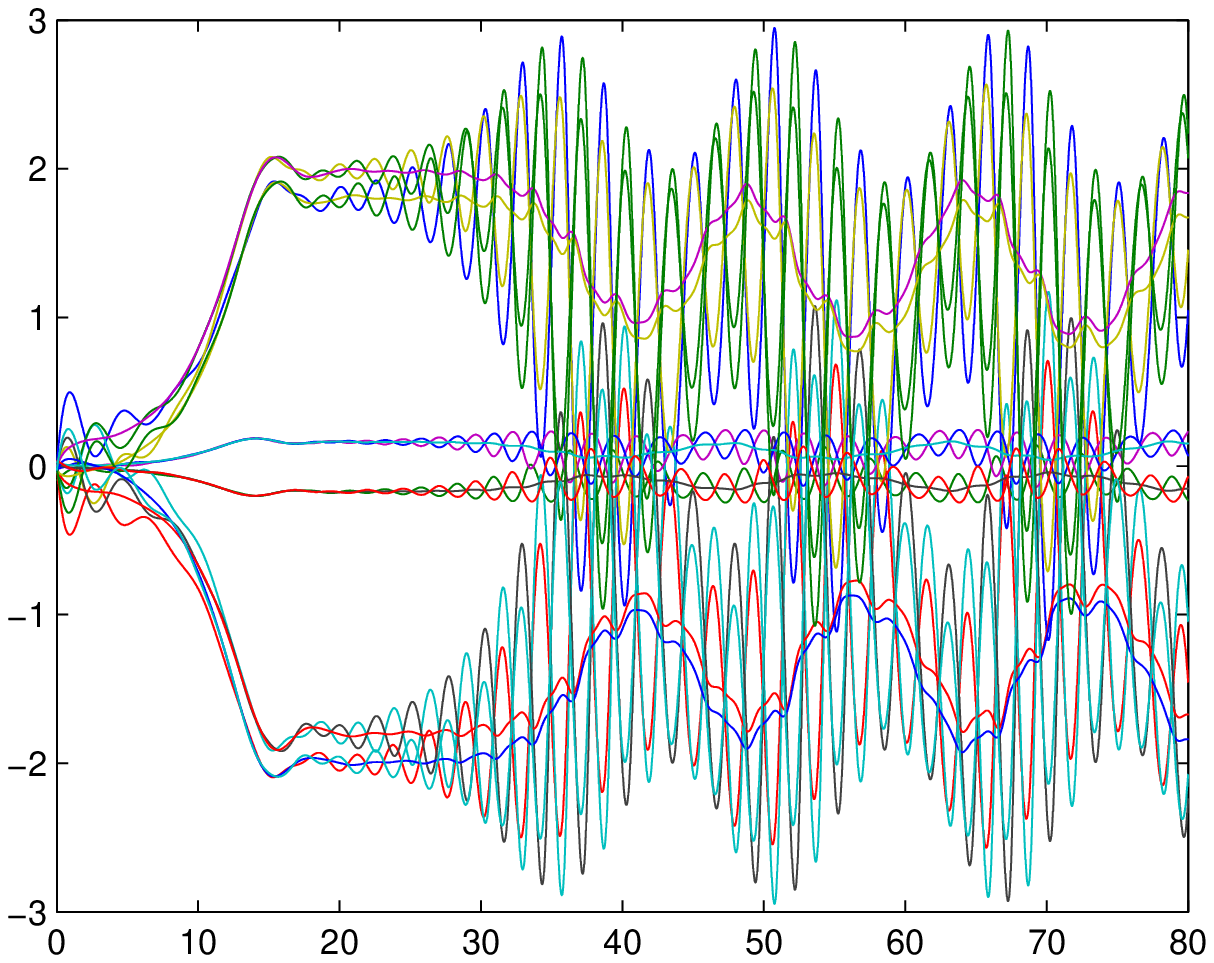}
 }
\end{center}
\vskip -0.3cm \centerline{ Fig. 5  Network on graph $G_1$. \qquad
\quad Fig. 6 Network on graph $G_2$.}

From Figs. 5 and  6, one can see that different topological
structures result in very different dynamical behaviors for the
smooth Chua's circuit network.

{\bf Example 2}\,\, Consider designing the inner coupling matrix
such that the corresponding network has an unbounded synchronized
region, where the node system is as in Example 1. In this case, $F$
is stable, so for any column vector $b$, $(F,b)$ is stabilizable,
e.g.,  $b=(0, 0, 1)^T$. By the method of Theorem 1, one gets
$k=(0.0708,   -0.15590,  0.4296).$ Then, change the inner coupling
matrix $H$ in Example 1 to $H=bk$, while keeping the other
parameters unchanged. This $F+\alpha H$ is stable for $\alpha \in
(-\infty, -1]$. Consequently, the corresponding network has an
unbounded synchronized region.  Figs. 7 and  8 show their similar
synchronization behaviors of the two networks built on graphs $G_1$
and $G_2$, respectively.

  \begin{center}
 \unitlength=1cm
 \qquad \hbox{\hspace*{0.1cm}  \epsfxsize5.5cm \epsfysize5cm
\epsffile{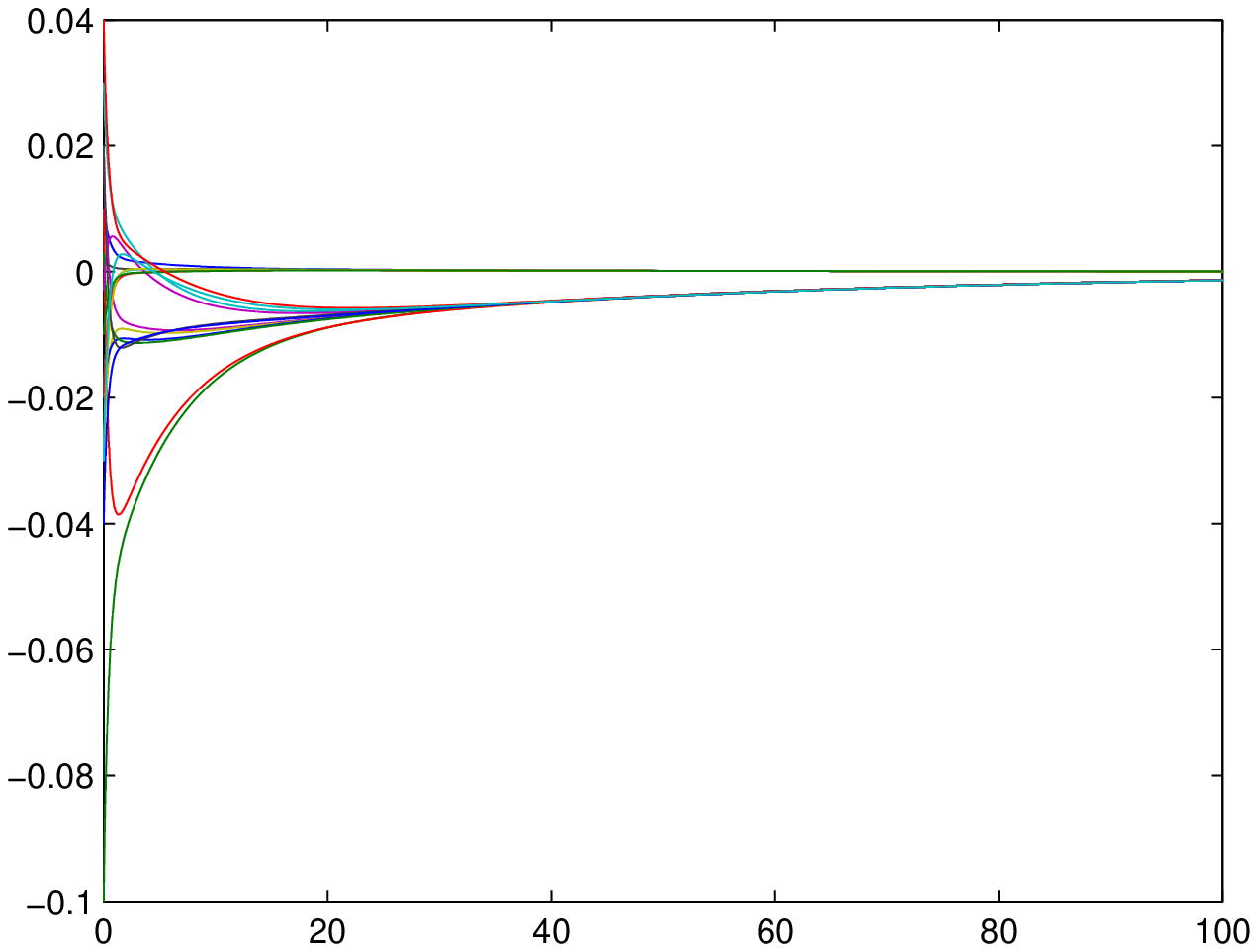}\qquad \epsfxsize5.5cm \epsfysize5cm
\epsffile{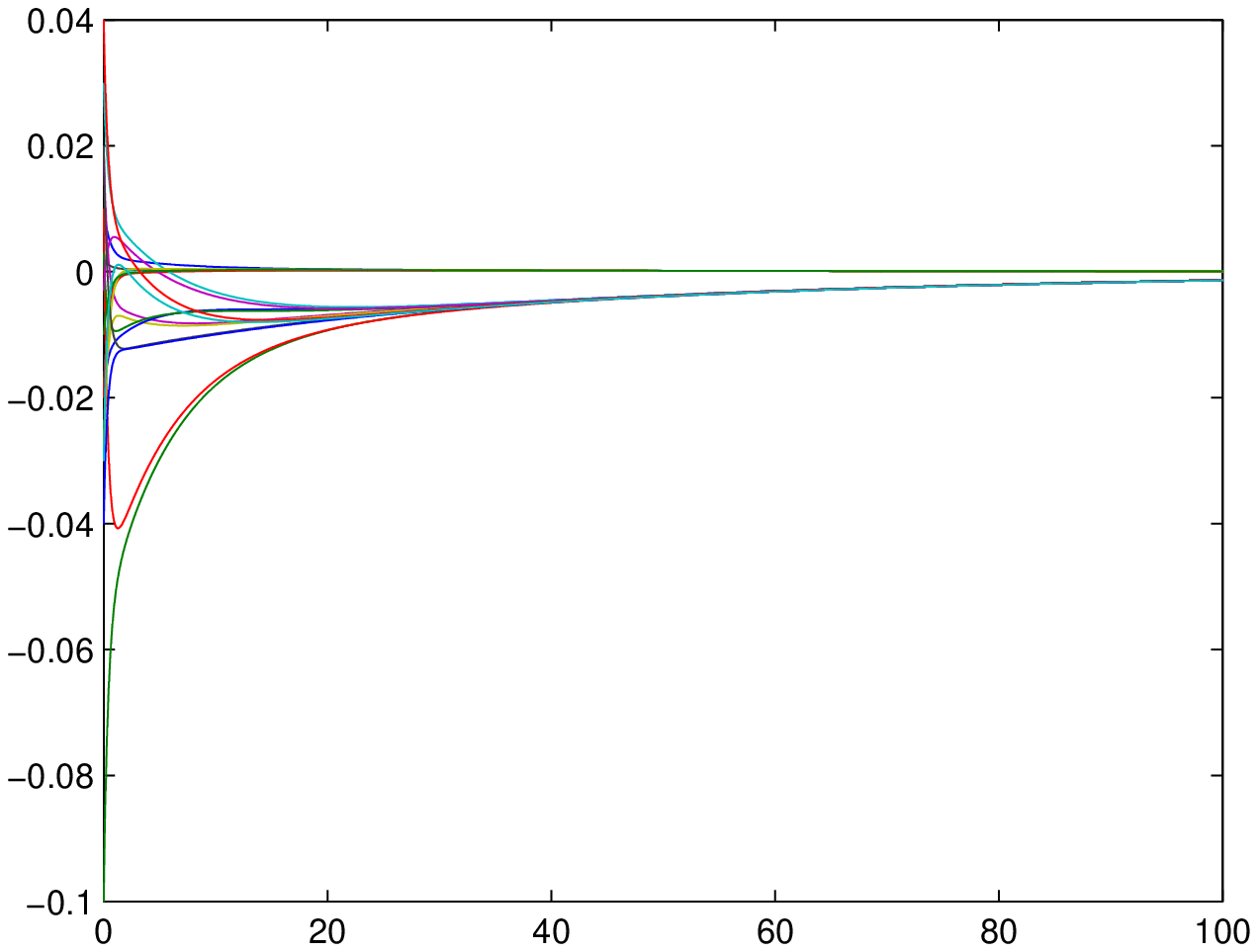}
 }
\end{center}
\vskip -0.3cm \centerline{ Fig. 7 Network on graph $G_1$. \qquad
\quad Fig. 8 Network on graph $G_2$.}

One can see that the network solutions synchronize a little faster
in Fig. 7 than that  in Fig. 8. This also demonstrates that the
synchronizability of the network built on $G_1$ is better than the
one  on $G_2$, which verifies the theoretical analysis given in the
previous sections.

\section{Conclusion}

In this paper, the synchronizability of complex dynamical  networks,
which is directly related to the inner linking matrix and the
topological matrix, has been carefully discussed. Two simple graphs
have been given to show that networks can have different
synchronizabilities even when they have the same average distance,
node betweenness centrality and degree distribution. It has also
been shown that the larger clustering coefficient, smaller
betweenness centrality and average distance  do not necessarily
imply better synchronizability. This demonstrates the complexity in
the relationship between the synchronizability and network
structural parameters. The most significant discovery of this paper
is that if the synchronized region is bounded, adding edges can
either increase or decrease the network synchronizability; however,
for networks with disconnected complementary graphs, adding edges
never decreases their synchronizability. Therefore,  better
understanding and careful manipulation of complementary graphs is
useful for enhancing network synchronizability. Moreover,  unbounded
synchronized regions are easier to analyze than the bounded ones. To
effectively enhance the network synchronizability,  a design method
for the inner linking matrix of rank 1 is finally provided such that
the resulting network has an unbounded synchronized region for the
case where the synchronous state is an equilibrium point of the
network.

 \hspace*{34pt}


\begin{thebibliography}{1}



\bibitem{bar02} M. Barahona, L. M. Pecora. Synchronization in small-world systems,
 {\it Phys. Rev. Lett.,} 89(5): 054101, 2002.

\bibitem{bel05} I. V. Belykh, E. Lange, M. Hasler. Synchronization of bursting neurons: what matters in
 the network topology, {\it Phys. Rev. Lett.,} 94: 188101,
2005.

\bibitem{boc06} S. Boccaletti, V. Latora, Y. Moreno, M. Chavez, D. U. Hwang.
Complex networks: structure and dynamics,
 {\it Physics Reports,} 424: 175-308, 2006.

\bibitem{cur97}  P. F. Curra, L. O. Chua. Absolute stability theory and
synchronization problem, {\it International Journal of Bifurcation
and Chaos}, 7: 1375-83, 1997.

\bibitem{duan07}  Z. S.  Duan, G. Chen, L. Huang.  Disconnected synchronized regions of
complex dynamical networks, {\it  submitted}, 2007.

\bibitem{duan04}  Z. S.  Duan, L. Huang, L. Wang,
J. Z. Wang.  Some applications of small gain theorem to
interconnected systems, {\it Syst. Contr. Lett.,}  52(3-4): 263-273,
2004.


\bibitem{hol02} P. Holme, B. J. Kim. Vertex overload breajdown in evolving networks.
{\it Phys. Rev. E,}  65: 066109, 2002.


\bibitem{hong04} H. Hong, B. J. Kim, M. Y. Choi, H. Park.
Factors that predict better synchronizability on complex networks,
{\it Phys. Rev. E,} 69: 067105, 2004.

\bibitem{iwa94} T. Iwasaki, R. E. Skelton.
All controllers for the general $H_\infty$ control problem: LMI
rxistence conditions and state space formulas, {\it Automatica,}
30(8): 1307-1317, 1994.

\bibitem{jia06} G. P. Jiang, W. K. S. Tang, G. R. Chen.
A state-observer-based approach for synchronization in complex
dynamcal networks. {\it IEEE Trans. Circuits  Syst.-I,} 53(12):
2739-2745, 2006.

\bibitem{koc05}
L. Kocarev, P. Amato. Synchronization in power-law networks, {\it
Chaos}, 15:  024101, 2005.

\bibitem{liu07} C. Liu, Z. S. Duan, G. Chen, L. Huang.
Analysis and control of synchronization regions in complex dynamical
networks, {\it submitted}, 2007.

\bibitem{luw04} W. Lu, T. Chen.
 Synchronization analysis of linearly coupled networks of discrete time systems,
 {\it Physica D,} 198: 148-168, 2004.

\bibitem{lu04}
J. H. L\"{u}, X. H. Yu, G.  Chen, D. Z. Cheng. Characterizing the
synchronizability of small-world dynamical networks, {\it IEEE.
Trans. Circuits Syst.-I }, 51(4), 787-796, 2004.

\bibitem{mer98} R. Merris. Laplacian graph eigenvectors, {\it Linear
Algebra and its Applications,} 278: 221-236, 1998.

\bibitem{mot05} A. E. Motter, C. S. Zhou, J. Kurths. Enhancing complex-network
synchronization, {\it Europhysics Letters,} 69(3): 334-340, 2005.

\bibitem{nis03}
T. Nishikawa, A. E. Motter, Y. C. Lai, F. C. Hoppensteadt.
Heterogeneity in oscillator netwroks: Are smaller worlds easier to
synchronize?  {\it Phys. Rev. Lett.,}  91(1): 014101, 2003.



\bibitem{pec98}
L. M. Pecora, T. L. Carroll, Master stability functions for
synchronized coupled systems, {\it Phys. Rev. Lett.}, 80(10):
2109-2112, 1998.




\bibitem{sor07} F. Sorrentino, M. di Bernardo, F. Garofalo, G. Chen.
 Controllability of complex networks via pinning,
 {\it Phys. Rev. E,} 75: 046103, 2007.

\bibitem{ste07} A. Stefa\'nski, P. Perlikowski, T. Kapitaniak.
 Ragged synchronizability of coupled oscillators, {\it Phys. Rev. E,}
75: 016210, 2007.

\bibitem{tsu05} A. Tsuneda. A gallery of attractors from smooth Chua's equation",  {\em
Int. J. Bifurcation and Chaos, }  15(1): 1-49, 2005.

\bibitem{wang06} X. Wang, Y. C. Lai, C. H. Lai.
 Synchronization is hierchical cluster networks,
 {\it arXiv. nlin. CD/0612057v1,}  2006.

\bibitem{wang02}  X. F. Wang, G. Chen.
 Synchronization in scale-free dynamical networks: robustness and fragility,
 {\it IEEE Trans.
Circuits Syst-I,} 49(1): 54-62, 2002.


\bibitem{wat98}  D. J. Watts, S. H. Strogatz.
 Collective dynamics of `small-world' networks,
 {\it Nature,} 393(6684): 440-442, 1998.

\bibitem{wu02}  C. W. Wu. Synchronization in coupled chaotic circuits
and systems. Singapore: World Scientific; 2002.

\bibitem{wu05} C. W. Wu. Synchronizability of networks of chaotic systems
coupled via a graph with a prescribed degree sequence, {\it Phys.
Lett. A,} 346: 281-287, 2005.

\bibitem{wu03} C. W. Wu. Perturbation of coupling matrices and its
effect on the synchronization in arrays of coupled chaotic systems,
 {\it Phys. Lett. A,} 319: 495-503, 2003.

\bibitem{zhao06} M. Zhao, B. H. Wang, G. Yan, H. J. Yang, W. J. Bai.
 Relations between average distance, heterogeneity and network synchronizability,
 {\it Physica A,} 371(2): 773-780, 2006.


\bibitem{zhou06} C. S. Zhou, J. Kurths. Dynamical weights and enhanced synchronization in adaptive complex
networks, {\it Phys. Rev. Lett.}, 96(16): 164102, 2006.

\end{thebibliography}
\end{document}